\documentclass[12pt]{article}
\usepackage[english]{babel}
\usepackage{amssymb,amsmath,amsthm}
\newtheorem{thm}{Theorem}[section]
\newtheorem{lem}{Lemma}[section]

\newtheorem{defn}{Definition}[section]

\begin{document}
\title{Classification of some nilpotent class\\ of Leibniz superalgebras}
\author{L.M. Camacho, J.R. G\'{o}mez, R.M. Navarro, B.A. Omirov}

\date{}
\maketitle

\begin{abstract}

The aim of this work is to present the description of Leibniz superalgebras up to
isomorphism with characteristic sequence
$(n | m-1, 1)$ and nilindex $n+m.$


\end{abstract}

\section{Introduction}

During many years the theory of Lie superalgebras has been
actively studied by many mathematicians and physicists. Many works
have been devoted to them but few have dealt with nilpotent Lie
superalgebras. Recent works \cite{Gilg1}--\cite{Gomez3} have
studied the problems of description of some classes of nilpotent
Lie superalgebras. It is well known that Lie superalgebras are a
generalization of Lie algebras \cite{Kac}. In the same way, the
notion of Leibniz algebras can be generalized by means of Leibniz
superalgebras. The elementary properties of Leibniz superalgebras
were obtained in \cite{O1}. The description of case of maximal
nilindex for nilpotent Leibniz superalgebras (nilpotent Leibniz
superalgebras distinguished by the feature of being
single-generated) is not difficult and was done in \cite{O1}.
However, the next stage is very problematic, it consists of
Leibniz superalgebras with dimensions even and odd parts equal to
n and m, respectively, and  of nilindex $n+m$. It should be noted
that such Lie superalgebras were classified in \cite{Gomez1}. Due
to the great difficulty in solving the problem, some restrictions
on the characteristic sequence are added. The experience of using
the characteristic sequence in Lie and Leibniz algebras (even in
Lie superalgebras) leads us to choose a restriction on this
invariant. Since graded non-commutative identity in non Lie
Leibniz superalgebras does not hold, we usually have to solve many
technical tasks when describing Leibniz superalgebras, \cite{03}.

In a similar way to Leibniz algebras and Lie superalgebras cases, it
is possible to define the notions of null-filiform and filiform
Leibniz superalgebras \cite{O2}, \cite{Gilg1}  as superalgebras with
characteristic sequences $(n | m)$ and $(n-1, 1 | m)$ respectively.
We should take into account that the superalgebras in \cite{Gomez2}
show that all null-filiform superalgebras have nilindex $n+m$(except for a unique
superalgebra of maximal nilindex) and there
exist filiform superalgebras which also have nilindex $n+m$ for some
$n,$ $m.$ In the present paper we investigate Leibniz
superalgebras with the characteristic sequence $C(L)=(n | m-1,1)$
and with nilindex equal to $n+m.$

In this paper all spaces and superalgebras are considered over the
complex number field.

\section{Preliminaries}

We recall the definition of Leibniz superalgebras.

\begin{defn} A $\mathbb{Z}_2$-graded vector space
$L=L_0\oplus L_1$ is called \emph{a Leibniz superalgebra} if it is
equipped with a product $[ - , - ]$ which satisfies the following
conditions:
$$
[L_\alpha,L_\beta]\subseteq L_{\alpha+\beta(mod2)}\ \mbox{ for all
} \ \alpha,\beta \in \mathbb{Z}_2,
$$
$$
[x,[y,z]]=[[x,y],z]-(-1)^{\alpha\beta}[[x,z],y] \ \mbox{ -- graded
Leibniz identity }
$$
for all $x\in L,$ $y\in L_\alpha,$ $z\in L_\beta,$ $\alpha,\beta
\in \mathbb{Z}_2.$
\end{defn}
Note that if in $L$ the identity $[x,y]=-(-1)^{\alpha\beta} [y,x]$
(where $x\in L_{\alpha},$ $y\in L_{\beta}$ )holds, then the graded Leibniz and graded Jacobi identities
coincide. Thus, Leibniz superalgebras are a generalization of Lie
superalgebras.

Let us anote an example of non Lie Leibniz superalgebras, which
generalize the construction of non Lie Leibniz algebras \cite{Lod}.

Let $A=A_0\oplus A_1$ be an associative superalgebra over a field
$F$ and $D : A \to A$ be a $F$-linear map satisfying the
condition:
$$
D(a(Db))=DaDb=D((Da)b)
$$
for all $a, b \in A.$ If in the vector space $A$ we define the new
product:
$$
\langle a,b\rangle_{D} := a(Db)-(-1)^{\alpha\beta} D(b)a
$$
for $a \in A_\alpha,$ $b \in A_\beta,$ $A$ becomes a Leibniz
superalgebra.

Let us introduce some notations
$$
\Re(L)=\left\{ R_x\ |\ x\in L\right\},
$$ $$
Leib^{n,m}=\left\{ L=L_0\oplus L_1\ | \ dimL_0=n,\ dimL_1=m\right\}.
$$

It is not difficult to see that the set $\Re(L)$ will be a Lie
superalgebra with the following multiplication:
$$
\langle R_a, R_b\rangle := R_aR_b-(-1)^{\alpha\beta} R_bR_a
$$
for all $R_a\in \Re(L)_\alpha,$ $R_b\in \Re(L)_\beta.$

Let $V=V_0\oplus V_1,$ $W=W_0\oplus W_1$ be two
$\mathbb{Z}_2$-graded spaces. We say that a linear map $f : V \to
W$ has degree $\alpha$ (denoted as $deg(f)=\alpha$), if
$f(V_\beta)\subseteq W_{\alpha +\beta}$ for all $\beta\in
\mathbb{Z}_2.$

\begin{defn} Let $L$ and $L'$ be Leibniz
superalgebras. A linear map $f : L\to L'$ is called \emph{a
homomorphism} of Leibniz superalgebras if

1.  $f$ preserves the
grading, i.e. $f(L_0)\subseteq L_0'$ and $f(L_1)\subseteq L_1'$
($deg(f)=0$);

2. $f([x,y])=[f(x),f(y)]$ for all $x,y\in L.$\\
Moreover, if  $f$ is one-to-one then it is called \emph{an
isomorphism} of Leibniz superalgebras $L$ and $L'.$
\end{defn}

For a given Leibniz superalgebra $L$ we define a descending central
sequence as follows:
$$
L^1=L,\quad L^{k+1}=[L^k,L^1], \ k \geq 1.
$$

\begin{defn} A Leibniz superalgebra $L$ is called
nilpotent, if there exists  $s\in\mathbb N$ such that $L^s=0.$ The
minimal number $s$ with this property is called index of
nilpotency (nilindex) of the superalgebra $L.$
\end{defn}
\begin{defn} The set $R(L)=\left\{ z\in L\ |\ [L,
z]=0\right\}$ is called the right annihilator of a superalgebra
$L$.\end{defn}

Using the Leibniz graded identity it is not difficult to see that
$R(L)$ is an ideal of the superalgebra $L$. Moreover, elements of
the form $[a,b]+(-1)^{\alpha\beta}[b,a]$ ($a\in L_\alpha, \ b\in
L_\beta$) belong to $R(L)$.

The description of Leibniz superalgebras of maximal nilindex is
represented in the following theorem.

\begin{thm} \emph{\cite{O1}.} \label{thm21} Let $L$ be an $n$-dimensional Leibniz superalgebra with maximal index of nilpotency.
Then $L$ is isomorphic to one of the following two non isomorphic
superalgebras:
$$
[e_i,e_1]=e_{i+1},\quad 1\le i\le n-1
$$ $$
\left\{ \begin{array}{ll} {[}e_i,e_1{]}=e_{i+1},& 1\le i\le n-1 \\[2mm] {[}e_i,e_2{]}=2e_{i+2},& 1\le i\le n-2\\
\end{array}\right.
$$
where the omitted products are zero.
\end{thm}

It should be noted that for the second superalgebra when $n+m$ is
even, we have $m=n$ and if $n+m$ is odd then $m=n+1.$ Moreover, it
is clear that the Leibniz superalgebra has the maximal nilindex if
and only if it is one-generated.

We define the characteristic sequence as in \cite{Gilg1}.

Let $L=L_0\oplus L_1$ be a nilpotent Leibniz superalgebra. For an arbitrary
element $x\in L_0,$ the operator of right multiplication $R_x$ is a nilpotent
endomorphism of the space $L_i,$ where $i\in \{0, 1\}.$ Let us denote by
$C_i(x)$ ($i\in \{0, 1\}$) the descending sequence of the dimensions of Jordan
blocks of the operator $R_x.$ Consider the lexicographical order on the set
$C_i(L_0)$.

\begin{defn} A sequence
$$
C(L)=\left( \left.\max\limits_{x\in L_0\setminus [L_0,L_0]}
C_0(x)\ \right|\ \max\limits_{\widetilde x\in L_0\setminus
[L_0,L_0]} C_1\left(\widetilde x\right) \right)
$$
is said to be the characteristic sequence of the Leibniz
superalgebra $L.$
\end{defn}
Similarly to \cite{Gilg2} (corollary 3.0.1) it can be proved that
the characteristic sequence is invariant under isomorphisms.




\section{Description of Leibniz superalgebras with characteristic
sequence $(n \ |\ m-1, 1)$ and nilindex $n+m$}

The following theorem gives us the location of the generators of the Leibniz
superalgebra from $Leib^{n,m}$ with characteristic sequence $(n \ |\ m-1, 1)$
and nilindex $n+m$.

\begin{thm} \label{thm31}
Let $L=L_0\oplus L_1$ be a Leibniz superalgebra from the variety
$Leib^{n,m}$ with characteristic sequence $(n \ |\ m-1, 1)$ and
nilindex $n+m.$ Then $L$ is two-generated and they belong to
$L_1.$\end{thm}
\begin{proof}
Since the nilindex of $L$ is $n+m$, then superalgebra $L$ is
two-generated. From the definition of characteristic sequence we can
conclude that there exists a basis $\{y_1, y_2, \dots, y_m\}$ of
$L_1$ such that the operator ${R_{x_1}}_{|L_1}$ in this basis has
one of the following forms:
$$\left( \begin{array}{cc}
  J_{n-1} & 0 \\
  0 & J_1 \\
\end{array}%
\right), \ \ \left(%
\begin{array}{cc}
  J_1 & 0 \\
  0 & J_{n-1} \\
\end{array}%
\right).$$ By a change of the basis elements $\{y_1, y_2, \dots,
y_m\},$ we can assume that the operator ${R_{x_1}}_{|L_1}$ has the
first form.

Since $C(L_0)=(n),$ then from \cite{O2} (example 1) we have that
$L_0$ is a zero-filiform Leibniz algebra. Without loss of generality
we may suppose that
$$
\left\{\begin{array}{ll}  [x_i,x_1]=x_{i+1}, & 1\le i\le n-1, \\[2mm]
[y_j,x_1]=y_{j+1}, & 1\le j\le m-2, \\[2mm]
[y_m,x_1]=[y_m,x_1]=0,
\end{array}\right.
$$
where $\{x_1,x_2, \dots, x_n, y_1, y_2, \dots, y_m\}$ is the basis
of the superalgebra $L.$

From these multiplications we deduce that $\{x_2, x_3,\dots, x_n\}$
lie in the right annihilator of the superalgebra $L$ and generators
can be selected as a linear combination of elements of the set
$\{x_1, y_1, y_m\}$.

It is clear that two generators can not lie in $L_0.$ Let us suppose
the opposite assertion to the assertion of the theorem, i.e. one
generator element lies in $L_0$ and the second lies in $L_1$. Then
we can choose as generators the elements $x_1, \ Ay_1+By_m.$

\textbf{Case 1.} Let us suppose $A=0.$ Then $x_1$ and $y_m$ can be chosen as
generators. Hence
$$L^2=\{x_2, x_3, \dots, x_n, y_1, y_2, \dots, y_{m-1}\}.$$

Introduce the notations
$$[y_1,y_m]=\sum\limits_{k=2}^{n}\beta_kx_k, \quad
[x_1,y_m]=\sum\limits_{s=1}^{m-1}\alpha_sy_s.$$

Consider the products
$$[x_i,[y_m,x_1]]=[[x_i,y_m],x_1] - [[x_i,x_1],y_m], \ 1 \leq i
\leq n.$$

On the other hand, $[x_i,[y_m,x_1]]=0.$

Therefore
$$\begin{array}{ll}
[x_i,y_m]=\sum\limits_{k=1}^{m-i}\alpha_ky_{k-1+i}, & 1 \leq i \leq \min\{n, m-1\}, \\[1mm]
[x_i,y_m]=0, & \min\{n, m-1\} < i \leq n.
\end{array}$$

Thus, we obtain $y_1\notin \langle[x_i,y_m]\rangle, \ 2 \leq i \leq
n$ and since $y_1\in L^2$, then $\alpha_1\neq0.$

Consider the products
$$[y_j,[y_m,x_1]]=[[y_j,y_m],x_1] - [[y_j,x_1],y_m], \ 1 \leq j\leq
m-2.$$

On the other hand we have $[y_j,[y_m,x_1]]=0, \ \ 1 \leq j\leq
m-2.$

Therefore
$$\begin{array}{ll}
[y_j,y_m]=\sum\limits_{s=2}^{n+1-j}\beta_{s}x_{s-1+j}, & 1 \leq j
\leq \min\{m-1, n-1\}, \\[1mm]
[y_j,y_m]=0, & \min\{n-1, m-1\} < j \leq m-1.
\end{array}$$

Consider the equalities
$$[x_1,[y_m,y_m]]=2[[x_1,y_m],y_m]=2[\sum\limits_{s=1}^{m-1}\alpha_sy_s,y_m]=\sum\limits_{s=1}^{m-1}\alpha_s
\sum\limits_{t=2}^{n+1-s}\beta_tx_{t-1+s}.$$

Since $[y_m,y_m] \in \langle x_2, x_3, \dots, x_n\rangle$ then
$[y_m,y_m]\in R(L)$ and hence $[x_1,[y_m,y_m]]=0$.

Therefore,
$$\sum\limits_{s=1}^{m-1}\alpha_s
\sum\limits_{t=2}^{n+1-s}\beta_tx_{t-1+s}=0.$$

Comparing the coefficient at the basic elements $x_i, \ 2 \leq i
\leq n,$ we obtain $\beta_i=0, \ 2 \leq i \leq n.$ And we have
$L^3=\{x_3, x_4, \dots, x_n, y_2, y_3, \dots, y_{m-1}\}$, but it
follows that the index of nilpotency of the superalgebra $L$ is
smaller than $n+m.$ Thus, we have contradiction with assumption
$A=0.$

\textbf{Case 2.} $A\neq 0.$ Then we can take as generators $x_1$ and
$y_1$. Hence, $L^2=\{x_2, x_3, \dots, x_n, y_2, y_3, \dots, y_m\}.$

Let us define the notations
$$[x_i,y_1]=\sum\limits_{k=2}^{m}\alpha_{i,k}y_k, \ 1 \leq i \leq
n, \quad [y_j,y_1]=\sum\limits_{s=2}^{n}\beta_{j,s}x_s, \ 1 \leq j
\leq m. \eqno(3.1)$$

Let us suppose $x_2\notin L^3.$ Then $L^3=\{x_3, x_4, \dots, x_n,
y_2, y_3, \dots, y_m\}$ and there exist some $i_0, j_0 \ (2 \leq
i_0, j_0)$ such that
$\alpha_{i_0,2}\alpha_{j_0,m}-\alpha_{j_0,2}\alpha_{i_0,m}\neq 0$
and $\alpha_{i_0,2}\alpha_{j_0,m}\neq0.$

Since $[x_{i_0},y_1]=\sum\limits_{k=2}^{m}\alpha_{i_0,k}y_k \in
R(L)$ and $\alpha_{i_0,2}\neq0$ then multiplying from the right side
by $x_1$ enough times, we obtain $\{y_3, y_4, \dots,
y_{m-1}\}\subseteq R(L)$.

Therefore
$$[x_{i_0,2},y_1]-\sum\limits_{k=3}^{m-1}\alpha_{i_0,k}y_k =
\alpha_{i_0,2}y_2+ \alpha_{i_0,m}y_m \in R(L), $$
$$[x_{j_0,2},y_1]-\sum\limits_{s=3}^{m-1}\alpha_{j_0,s}y_s
=\alpha_{j_0,2}y_2+ \alpha_{j_0,m}y_m \in R(L).$$

Since $\alpha_{i_0,2}\alpha_{j_0,m}-\alpha_{j_0,2}\alpha_{i_0,m}\neq
0,$ then we have $y_2, y_m \in R(L),$ i.e. $\{y_2, y_3, \dots,
y_m\}\subseteq R(L)$ and hence $[x_1,y_1]\in R(L).$

From equality
$$[[x_{i-1},x_1],y_1]=[x_{i-1},[x_1,y_1]]+[[x_{i-1},y_1],x_1],$$
we obtain
$$\begin{array}{ll}
[x_i,y_1]=\sum\limits_{k=2}^{m-i}\alpha_{1,k}y_{k+i-1}, & 2 \leq
i \leq \min\{n, m-2\}, \\[1mm]
[x_i,y_1]=0, & \min\{n, m-2\} < i \leq n,
\end{array}$$ \\[1mm]
from which we have a contradiction to the assumption
$\alpha_{i_0,2}\neq0.$

Therefore $x_2 \in L^3$ and $x_2$ is a linear combination of
products $[y_i,y_1], \ 1 \leq i \leq m,$ hence
 $[x_2,y_1]=\sum\limits_{i=1}^m \gamma_i[[y_i,y_1],y_1].$ Using
the Leibniz graded identity and $[y_1,y_1]\in R(L)$ one can easy see that
$[x_2,y_1]=0.$

Thus, we have $L^3=\{x_2, x_3, \dots, x_n, Cy_2+Dy_m, y_3, \dots,
y_{m-1}\}.$

Since $[x_2,y_1]=0$, then $Cy_2+Dy_m \in L^3$ should be expressed by
linear combinations of products $[x_i,y_1], \ 3 \leq i \leq n$.
Therefore, $Cy_2+Dy_m \in L^5.$

Moreover, if  $C\neq 0,$ then $Cy_2+Dy_m \in R(L).$ Then multiplying
from the  right side by $x_1$ enough times, we obtain $\{y_3, y_4,
\dots, y_{m-1}\}\subseteq R(L)$.

Since $[x_2, y_{1}]=0$ and $Cy_2+Dy_m\in L^3$ then there exist
${t_0}, \ t_0\geq 3$, such that $[x_{t_0},
y_1]=c_2(Ay_2+By_m)+c_3y_3+\dots+c_{m-1}y_{m-1}, \ c_2 \neq 0.$

Applying the Leibniz graded identity we obtain

$$[x_1,[\dots [[y_1, \underbrace{x_1],x_1],\dots,
x_1}\limits_{t_0-1-times}]]=(-1)^{t_0-1}c_2(Cy_2+Dy_m)+\{y_3,
y_4,\dots, y_{m-1}\}. \eqno(3.2)$$

As $L$ is nilpotent, we have existence $s \in \mathbb{N},$ such that
$y_{t_0}\in L^s\setminus L^{s+1}.$

If $t_0<m,$ then $[x_1,y_{t_0}]=(-1)^{t_0-1}c_2(Cy_2+Dy_m)+\{y_3, y_4,\dots,
y_{m-1}\}\in L^{s+1}$ and multiplying this equality $(t_0-2)$ times from the
right side by $x_1,$ we obtain that $(-1)^{i_0-1}c_2Cy_{t_0} \in L^{s+t_0-1}.$
The inequality $s+t_0-1>s$ contradicts the condition $y_{t_0}\in L^s\setminus
L^{s+1}$ and hence we obtain $C=0.$

If $t_0\geq m,$ then $[x_1,y_{t_0}]=0.$ From (3.2) we again obtain
$C=0.$

Thus, $L^3=\{x_2, x_3,\dots, x_n, y_3, \dots, y_{m-1}, y_m \}.$

Consider the following subcases.

\textbf{Case 2.1} Let $\alpha_{1,2}\neq 0.$ Suppose that $x_2\in
L^l\setminus L^{l+1}$ for some $l \ (3\leq l \leq m).$ Then
$$L^l=\{x_2, x_3,\dots, x_n, y_l, \dots, y_{m-1}, y_m \}.$$
$$L^{l+1}=\{x_3, x_4,\dots, x_n, y_l, \dots, y_{m-1}, y_m
\}.$$

Since $x^2 \in L^l\setminus L^{l+1}$, we obtain $\beta_{l-1,2}\neq
0,$ i.e.
$$[y_{l-1},y_1]=\sum\limits_{k=2}^n\beta_{l-1,k}x_k.$$

The equality $[x_2,y_1]=0,$ deduce $y_l \in L^{l+2}$ and
$L^{l+2}=\{x_4, x_5,\dots, x_n, y_l, \dots, y_{m-1}, y_m \}.$

Therefore $[y_l,y_1]=\sum\limits_{k=4}^n\beta_{l,k}x_k,$ i.e.
$\beta_{l,2} = \beta_{l,3} = 0$

In notation (3.1) by induction on $j$ for any value of $i$ one can prove the
following equality:
$$
[y_i,y_j]=\sum\limits_{k=0}^{\min\{ i+j-1,m-1\}-i} (-1)^k
C^k_{j-1} \sum\limits_{t=2}^{n-j+k+1} \beta_{i+k,t} x_{t+ j-k-1},
\eqno (3.3)
$$
where $1 \le j \le m-1,$ $1 \le i \le m-1.$

From (3.3) we have $[y_2,y_l]=\beta_{l-1,2}x_3+{x_4, x_5,\dots,
x_n}$

Consider the equalities
$$[x_1,[y_1,y_l]]=[[x_1,y_1],y_l]+[[x_1,y_l],y_1]=$$
$$=\sum\limits_{k=2}^{m}\alpha_{1,k}[y_k,y_l]+
\sum\limits_{s=l+1}^{m}\gamma_{1,s}[y_k,y_1]=\alpha_{1,2}\beta_{l-1,2}x_3+\{x_4,
..., x_n\}.$$

On the other hand $$[x_1,[y_1,y_l]]=0.$$

So $\alpha_{1,2}\beta_{l-1,2}x_3=0,$ i.e. we have a contradiction
with supposition $\alpha_{1,2} \neq 0.$

\textbf{Case 2.1.} If $\alpha_{1,2}=0.$ Then
$[y_1,x_1]+[x_1,y_1]=y_2+\alpha_{1,3}y_3+\alpha_{1,4}y_4+\dots+\alpha_{1,m}y_m
\in R(L)$ and then multiplying from the right side by $x_1$ enough times, we
obtain $y_2+\alpha_{1,m}y_m, y_3,\dots, y_{m-1}\in R(L).$

Since $y_m \in L^3$, there exists $i_0\geq 2$ such that

$[x_{i_0},y_1]=\alpha_{{i_0},3}y_3+\alpha_{{i_0},4}y_4+\dots+\alpha_{{i_0},m}y_m,$
$\alpha_{{i_0},m} \neq 0.$

From $[x_{i_0},y_1] \in R(L),$ we have
$\alpha_{{i_0},3}y_3+\alpha_{{i_0},4}y_4+\dots+\alpha_{{i_0},m}y_m,
\in R(L).$ Therefore $y_2, y_m \in R(L).$

Consider the equalities

$[y_j,[y_1,x_1]]=[[y_j,y_1],x_1]-[[y_j,x_1],y_1]=
[\sum\limits_{k=2}^{n+1-j}
\beta_{j,k}x_{k-1+j},x_1]-[y_{j+1},y_1]=$

$=\sum\limits_{k=2}^{n-j} \beta_{j,k}x_{k+j}-[y_{j+1},y_1].$

On the other hand we have $[y_j,[y_1,x_1]]=[y_j,y_2]=0.$

Therefore

$[y_{j+1},y_1]=\sum\limits_{k=2}^{n-j} \beta_{1,2}x_{k+j}, \ 1\leq
j \leq \min\{n-2,m-2\},$

$[y_{j+1},y_1]=0 \ \min\{n-2,m-2\}\leq j \leq m-2.$

From these we conclude that $x_2$ cannot be expressed via linear
combination of $\{[y_j,y_1]\}, 2\leq j\leq m.$ Therefore $x_2\notin
L^3,$ it is a contradiction with assumption $x_2\in L^3.$

Thus, the supposition of the case when one generator lies in $L_0$
and the second lies in $L_1$ leads to a contradiction, therefore
both generators belong to $L_1$ and the proof of the theorem is
completed.
\end{proof}

The following lemma is an important technical result in our
description.
\begin{lem}   \label{lem31}
Let $L=L_0\oplus L_1$ be a Leibniz superalgebra from the variety
$Leib^{n,m}$ with characteristic sequence $(n \ |\ m-1, 1)$ and
nilindex $n+m.$ Then there exists $y\in L_1$ such that the
superalgebra $\langle y \rangle$ is isomorphic to the superalgebra
from Theorem \ref{thm21} and $m\in \{ n+1, n+2\}.$
\end{lem}
\begin{proof}
From Theorem \ref{thm31} we have that both generators are in $L_1$. Let $\{
x_1,x_2,\ldots,x_n,y_1,y_2,\ldots,y_m\}$ be the basis of $L$ such that $y_1,$
$y_m$ are the generators and
$$[x_i,x_1]=x_{i+1}, \ 1 \leq i \leq n-1, \ [y_j,x_1]=y_{j+1}, \ 1
\leq j \leq m-2, $$
$$[y_1,y_1]=\sum\limits_{i=1}^n a_ix_i, \
 [y_m,y_1]=\sum\limits_{i=1}^n b_ix_i,$$
$$[y_1,y_m]=\sum\limits_{i=1}^n c_ix_i, \
[y_m,y_m]=\sum\limits_{i=1}^n d_ix_i.$$

Since $[y_j,x_1]=y_{j+1}, \ 1 \leq j \leq m-2, $ then $x_1$ cannot be generated
via multiplications of the elements $\{y_2, y_3, \dots, y_{m-1}\}.$ Therefore
$x_1\in \langle [y_1,y_1], [y_m,y_1], [y_1,y_m], [y_m,y_m] \rangle$ and hence
$(a_1, b_1, c_1, d_1\neq (0, 0, 0, 0).$

Let us suppose that $a_1\neq0.$ Then making the change
$$x_i^{'}=\sum\limits_{k=1}^{n+1-i}a_kx_{k+i-1}, \ 1 \leq i \leq
n, \ y_j^{'}=a_1^{i-1}y_j, \ 1 \leq j \leq m-1$$ we obtain
$$[x_i,x_1]=x_{i+1}, \ 1 \leq i \leq n-1, \ [y_j,x_1]=y_{j+1}, \ 1
\leq j \leq m-2, \ [y_1,y_1]=x_1.$$

It is clear that the subsuperalgebra $\langle y_1\rangle=\{x_1, x_2,
\dots, x_n, y_1, y_2, \dots, y_{m-1}\}$ is one-generated and
therefore it has maximal index of nilpotency. Hence from Theorem
\ref{thm21} we have that either $m=n+1$ or $m=n+2.$

Let us suppose now that $a_1=0.$

Consider the product
$$[y_m,[y_m,x_1]]=[[y_m,y_m],x_1] - [[y_m,x_1],y_m].$$
Then $[y_m,y_m]=d_nx_n$ and $d_i=0, \ 1 \leq i \leq n-1.$ It leads
to $(b_1,c_1)\neq(0, 0).$

From $[y_m,y_1]-[y_1,y_m]=\sum\limits_{i=1}^{n}(b_i-c_i)x_i \in
R(L)$ and $x_j\in R(L), \ 2 \leq j \leq n$ we obtain that
$(b_1-c_1)x_1\in R(L),$ but $x_1\notin R(L)$ and hence
$b_1=c_1\neq0.$

By the following change of basis:

$$x_1^{'}=\sum\limits_{k=1}^nb_kx_k, \
x_{i+1}^{'}=[x_i^{'},x_1^{'}], \ 1 \leq i \leq n-1, y_j^{'}=y_j, \
1 \leq j \leq  m$$

we can assume $[y_m,y_1]=x_1.$

From the products

$$[y_m,[y_i,x_1]]=[[y_m,y_i],x_1] - [[y_m,x_1],y_i],$$
$$[y_m,[y_1,y_1]]=2[[y_m,y_1],y_1]=2[x_1,y_1]$$
we obtain $[y_m,y_{i+1}]=x_{i+1}, \ 1 \leq i \leq m-2, \
[x_1,y_1]=0.$

Since $[y_1,x_1]+[x_1,y_1]\in R(L)$ then $y_2\in R(L),$ but this
contradicts $[y_m,y_2]=x_2$ and therefore the case $a_1=0$ is not
possible either.
\end{proof}

The existence of an adapted basis with conditions of Lemma \ref{lem31} in case
of $m=n+1$ is described in the following lemma.
\begin{lem}   \label{lem32}
Let $L$ be a Leibniz superalgebra from the variety $Leib^{n,m}$ with
characteristic sequence $(n \ |\ m-1, 1)$ and nilindex $n+m.$ Then,
in case $m=n+1,$ there exists a basis $\{ x_1, x_2, \ldots, x_n,
y_1, y_2, \ldots, y_{n+1}\}$ of $L$ in which the products have
the following form:\\[1mm]

$\begin{array}{lll} [x_i,x_1]=x_{i+1},\ 1\le i\le n-1,  &
[y_j,x_1]=y_{j+1},\ 1\le j\le n-1, \\[1mm]
[x_i,y_1]=\frac12 y_{i+1},\ 1\le i\le n-1,  & [y_j,y_1]=x_{j},\ 1\le j\le n, \\[1mm]
[y_{n+1},y_{n+1}]=\gamma x_n, &
[x_i,y_{n+1}]=\sum\limits_{k=\left[\frac{n+4}2\right]}^{n+1-i}
\beta_k y_{k-1+i},\ 1\le i\le \left[ \frac{n-1}2\right], \\[1mm]

[y_1,y_{n+1}]=-2\sum\limits_{k=\left[\frac{n+4}2\right]}^{n}
\beta_k x_{k-1}+\beta x_n,  &
[y_j,y_{n+1}]=-2\sum\limits_{k=\left[\frac{n+4}2\right]}^{n+2-j}
\beta_k x_{k-2+j},\ 2\le j\le \left[\frac{n+1}2\right]. \\[1mm]
\end{array}$

\end{lem}

\begin{proof}
From Lemma \ref{lem31} and Theorem \ref{thm21} we have that there
exists a  basis $\{ x_1, x_2, \ldots, x_n, y_1, y_2, \ldots,
y_{n+1}\}$ in which products have the following form:\\[1mm]

$\begin{array}{ll}
[x_i,x_1]=x_{i+1},\ 1\le i\le n-1, & [y_j,x_1]=y_{j+1},\ 1\le j\le n-1, \\[1mm]
[x_i,y_1]=\frac12 y_{i+1},\ 1\le i\le n-1, & [y_j,y_1]=x_{j},\ 1\le j\le n,\\[1mm]
[y_1,y_{n+1}]=\sum\limits_{k=1}^n\alpha_kx_k, & [x_1,y_{n+1}]=\sum\limits_{s=2}^n\beta_sy_s.\\[1mm]
\end{array}$

From equalities:
$$
[y_{n+1},[y_1,x_1]]=[[y_{n+1},y_1],x_1] - [[y_{n+1},x_1],y_1],
$$ $$
[y_{n+1},[x_1,y_1]]=[[y_{n+1},x_1],y_1] - [[y_{n+1},y_1],x_1],
$$
$$
[y_{n+1},y_{i+1}]=[y_{n+1},[y_i,x_1]]=[[y_{n+1},y_i],x_1] -
[[y_{n+1},x_1],y_{i}], \ 1 \leq i \leq n-1,
$$
we have $[y_{n+1},y_2]=[[y_{n+1},y_1],x_1]$ and
$[y_{n+1},y_2]=-2[[y_{n+1},y_1],x_1]$, it follows that

$[y_{n+1},y_1]=bx_n$ and $[y_{n+1},y_i]=0$ for $2\le i\le n.$

Consider the equality:
$$
[y_{n+1},[y_{n+1},x_1]]=[[y_{n+1},y_{n+1}],x_1]+
[[y_{n+1},x_1],y_{n+1}],
$$
then $[y_{n+1},y_{n+1}]=\gamma x_n.$

Considering the graded identities:
$$
[y_j,[y_{n+1},x_1]]=[[y_j,y_{n+1}],x_1] - [[y_j,x_1],y_{n+1}],
$$ $$
[x_i,[y_{n+1}],x_1]=[[x_i,y_{n+1}],x_1] - [[x_i,x_1],y_{n+1}],
$$
we have
$$[y_j,y_{n+1}]=\sum\limits_{k=1}^{n+1-j}\alpha_kx_{k-1+j}, \ 2
\leq j \leq n,$$
$$[x_i,y_{n+1}]=\sum\limits_{s=2}^{n+1-i}\beta_sy_{s-1+i}, \
2\leq i \leq n-1.$$ From the chain of equalities
$$0=[y_1,[y_{n+1},y_1]]=[[y_1,y_{n+1}],y_1] +
[[y_1,y_1],y_{n+1}]=$$ $$=[\sum\limits_{k=1}^{n}\alpha_kx_k,y_1] +
[x_1,y_{n+1}]= \frac{1}{2}\sum\limits_{k=1}^{n-1}\alpha_ky_{k+1} +
\sum\limits_{s=2}^{n}\beta_sy_s,$$
we obtain $\alpha_i=-2\beta_{i+1}, \ 1 \leq i \leq n-1.$

Substituting these relations in the chain of equalities
$$0=[y_1,[y_{n+1},y_{n+1}]]=2[[y_1,y_{n+1}],y_{n+1}],$$
 we get $\beta_i=0, \ 2 \leq i \leq [\frac{n+2}{2}].$
Thus, we obtain the multiplications of the superalgebra as in the assertion of the lemma. \end{proof}

For convenience, we will denote the superalgebra from the family of
Lemma \ref{lem32} as $L\left( \gamma, \beta_{\left[
\frac{n+4}2\right]}, \beta_{\left[ \frac{n+4}2\right]+1}, \ldots,
\beta_n,\beta\right).$

Using the properties of adapted basis we obtain necessary and sufficient
conditions when two arbitrary superalgebras from the family of Lemma
\ref{lem32} are isomorphic.

\begin{thm} \label{thm32} Two superalgebras $L\left( \gamma, \beta_{\left[ \frac{n+4}2\right]}, \beta_{\left[ \frac{n+4}2\right]+1},
\ldots, \beta_n,\beta\right)$ and \\ $L'\left( \gamma',
\beta'_{\left[ \frac{n+4}2\right]}, \beta'_{\left[
\frac{n+4}2\right]+1}, \ldots, \beta'_n,\beta'\right)$ are
isomorphic if and only if there exist $a_1, a_{n+1}, b_{n+1}\in \mathbb{C}$
such that the following conditions hold:

for odd $n$:
$$
\left\{ \begin{array}{l} b^2_{n+1}\gamma=\gamma'a_1^{2n},\\[2mm] b_{n+1}\beta_j=a_1^{2j-3}\beta_j',\quad \left[
\frac{n+4}2\right]\le j\le n, \\[2mm] a_{n+1}b_{n+1}\gamma + a_1b_{n+1}\beta = a_1^{2n} \beta' + 4
\beta'_{\left[\frac{n+4}2\right]}a_1^{2 \left[\frac{n+1}2\right] -1}a_{n+1} \beta_{\left[\frac{n+4}2\right]},\\
\end{array}\right.
$$

for even $n$:
$$
\left\{ \begin{array}{l} b^2_{n+1}\gamma=\gamma'a_1^{2n},\\[2mm] b_{n+1}\beta_j=a_1^{2j-3}\beta_j',\quad \left[
\frac{n+4}2\right]\le j\le n, \\[2mm] a_{n+1}b_{n+1}\gamma + a_1b_{n+1}\beta = a_1^{2n} \beta' .\\
\end{array}\right.
$$
\end{thm}
\begin{proof}
Let us make a general change of generator elements in the form:
$$
y'_1=\sum\limits_{i=1}^{n+1}a_iy_i, \
y'_{n+1}=\sum\limits_{j=1}^{n+1}b_jy_j,
$$
where the rank $\left(\begin{array}{llll} a_1 & a_2 & \ldots & a_{n+1} \\ b_1 & b_2 & \ldots & b_{n+1} \\
\end{array}\right)=2.$

We express the new basis $\{ x_1',x_2',\ldots,x_n',y_1',
y_2',\ldots,y_{n+1}'\}$  of the superalgebra $L'\left( \gamma',
\beta'_{\left[ \frac{n+4}2\right]}, \beta'_{\left[
\frac{n+4}2\right]+1}, \ldots, \beta'_n,\beta'\right)$ with
respect to the old basis $\{ x_1, x_2, \ldots, x_n, y_1, y_2, \\
\ldots, y_{n+1}\}.$

Then for element $x_1'$ we have
$$
x_1'= [y_1',y_1']=a_1\sum\limits_{i=1}^n
a_ix_i-2a_1a_{n+1}\sum\limits_{k=\left[\frac{n+4}2\right]}^n\beta_kx_{k-1}
+$$ $$+a_1 a_{n+1} \beta x_n- 2a_{n+1}
\sum\limits_{i=2}^{\left[\frac{n+1}2\right]}a_i
\sum\limits_{k=\left[\frac{n+4}2\right]}^{n +2-i}\beta_kx_{k+i-2}+
a^2_{n+1}\gamma x_n.
$$

The expression of $x_{t+1}'$ $\left( 1\le t\le
\left[\frac{n-1}2\right]\right)$ will be as follows:
$$
x_{t+1}'=\left[ x_{t}', x_1'\right]
=a_1^{2t+1}\sum\limits_{i=1}^{n-t}a_ix_{i+t}-2a_1^{2t}a_{n+1}
\sum\limits_{i=1}^{ \left[\frac{n+1}2 \right]-t} a_i
\sum\limits_{k=\left[\frac{n+4}2\right]}^{n+2-t-i}\beta_k
x_{k+t+i-2}.
$$

And for $x_{t+1}'$ $\left( \left[\frac{n+1}2\right]\le t\le n-1 \right)$ we have
$$
x_{t+1}'=\left[ x_{t}', x_1'\right] =a_1^{2t+1}\sum\limits_{i=1}^{n-t}a_ix_{t+i}.
$$

For basis elements $y_i'$ of the space $L_1'$ basis we have:
$$
y_t'=\left[ y_{t-1}',x_1'\right]= a_1^{2(t-1)} \sum\limits_{j=
1}^{n+1-t}a_iy_{t+j-1},\ 2\le t\le n.
$$

Consider the equalities:
$$
\left[ y_{n+1}',y_{n+1}'\right]=b_1\sum\limits_{i=1}^n
b_ix_i-2b_1b_{n+1}\sum\limits_{k=\left[\frac{n+4}2\right]}^n
\beta_kx_{k-1}+ b_1b_{n+1} \beta x_n-$$
$$
-2b_{n+1}\sum\limits_{i=2}^{\left[\frac{n+1}2\right]}b_i
\sum\limits_{k=\left[\frac{n+4}2\right]}^{n +2-i}\beta_kx_{k+i-2}
+ b_{n+1}^2\gamma x_n=\gamma'x_n'=\gamma'a_1^{2n}x_n. $$

Comparing the coefficients we obtain the following restrictions:
$$
b_1=b_2=\ldots=b_{\left[\frac{n-1}2\right]}=0, \eqno (3.4)
$$ $$
-2b_{n+1}
b_{\left[\frac{n+1}2\right]}\beta_{\left[\frac{n+4}2\right]}
+b_{n+1}^2\gamma =\gamma' a_1^{2n}. \eqno (3.5)
$$

Using the chain of the equalities:
$$
\left[ y_{n+1}',y_1'\right]=
a_1\sum\limits_{i=\left[\frac{n+1}2\right]}^nb_ix_i+a_{n+1}b_{n+1}\gamma
x_n- 2a_{n+1}
b_{\left[\frac{n+1}2\right]}\beta_{\left[\frac{n+4}2\right]}
x_n=0,
$$
the equalities (3.4) - (3.5) and comparing the coefficients to the
basis elements we obtain:
$$
\left\{\begin{array}{l} b_1=b_2= b_{\left[\frac{n-1}2\right]} = b_{\left[\frac{n+1}2\right]}= \ldots=b_{n-1}
=0,\\[2mm] b^2_{n+1}\gamma=\gamma'a_1^{2n}, \\ a_1b_n+a_{n+1}b_{n+1}\gamma=0.\\[2mm]
\end{array} \right.\eqno (3.6)
$$

Therefore, $y'_{n+1}=b_ny_n+b_{n+1}y_{n+1}.$

Consider the products:
$$
\left[ x_{1}',y_{n+1}'\right]= \left[ a_1\sum\limits_{i=1}^na_ix_i-2a_1a_{n+1} \sum\limits_{k=\left[\frac{n+4}2
\right]}^n \beta_kx_{k-1}+a_1a_{n+1}\beta x_n - \right.
$$ $$
\left.-2a_{n+1}\sum\limits_{i=2}^{\left[\frac{n+1}2\right]} a_i
\sum\limits_{k=\left[\frac{n+4}2\right]}^{n+2-i}\beta_k x_{k+i
-2}+a_{n+1}^2 \gamma x_n, b_ny_n+b_{n+1}y_{n+1}\right] =
$$

$ =a_1b_{n+1} \sum\limits_{i=1}^{\left[\frac{n-1}2\right]}a_i
\sum\limits_{k=\left[\frac{n+4}2\right]}^{n+1-i} \beta_k
y_{k+i-1}, $
$$ \left[ x_{1}',y_{n+1}'\right]=
\sum\limits_{k=\left[\frac{n+4}2\right]}^{n} \beta_k' y_{k}' =
\beta'_{\left[ \frac{n+4}2\right]} a_1^{2\left[
\frac{n+4}2\right]-2} \left( a_1 y_{\left[ \frac{n+4}2\right]}
+a_2 y_{\left[ \frac{n+4}2\right]+1} +\ldots + a_{n+1-\left[
\frac{n+4}2\right]}y_n\right)+
$$ $$
+\beta'_{\left[ \frac{n+4}2\right]+1} a_1^{2\left[
\frac{n+4}2\right]} \left( a_1 y_{\left[ \frac{n+4}2\right]+1}
+a_2 y_{\left[ \frac{n+4}2\right]+2} +\ldots + a_{n-\left[
\frac{n+4}2\right]}y_n \right)+ \ldots+
$$ $$
+\beta_n'a_1^{2n-1}y_n= \sum\limits_{\left[
\frac{n+4}2\right]}^n\beta_k'a_1^{2(k-1)}
\sum\limits_{j=1}^{n+1-k} a_i y_{k+j-1}.
$$

From which, comparing the coefficients we have the following
restrictions:
$$
b_{n+1}\beta_j=a_1^{2j-3}\beta_j', \qquad \left[
\frac{n+4}2\right]\le j \le n.
$$

Consider the following product on the one hand:
$$
\left[ y_1',y_{n+1}'\right]= [a_1y_1+a_2y_2+\ldots+a_{n+1}y_{n+1}, b_ny_n+b_{n+1}y_{n+1}]=
$$ $$
=-2a_1b_{n+1} \sum\limits_{k=\left[ \frac{n+4}2\right]}^n \beta_k
x_{k-1} +a_1b_{n+1}\beta x_n -2b_{n+1} \sum\limits_{i=2}^{\left[
\frac{n+1}2\right]} a_i \sum\limits_{k=\left[
\frac{n+4}2\right]}^{n+2-i} \beta_k x_{k+i-2}
+a_{n+1}b_{n+1}\gamma x_n
$$
and on the other hand, let us consider the followings product in the
case of an odd $n$.
$$
\left[ y_1',y_{n+1}'\right]= -2\sum\limits_{k=\left[
\frac{n+4}2\right]}^n \beta_k' x_{k-1}'+\beta'x_n'=
-2\beta_{\left[ \frac{n+4}2\right]}'\left( a_1^{2\left(\left[
\frac{n+4}2\right]-1\right)} x_{\left[ \frac{n+4}2\right] -1}+
\right.
$$ $$
+ \left. a_1^{2\left(\left[ \frac{n+4}2\right]-2\right)+1} a_2
x_{\left[ \frac{n+4}2\right]} + \ldots + a_1^{2\left( \left[
\frac{n+4}2\right]-2\right)+1} a_{n-\left[ \frac{n+4}2\right]+2}
x_n -2a_1^{2\left[ \frac{n-1}2\right]+1} a_{n+1} \beta_{\left[
\frac{n+4}2\right]}x_n \right)-$$

$$
-2\beta_{\left[ \frac{n+4}2\right]+1}'\left( a_1^{2\left[
\frac{n+4}2\right]} x_{\left[ \frac{n+4}2\right]}+ \right. \left.
a_1^{2\left(\left[ \frac{n+4}2\right]-1\right)+1} a_2 x_{\left[
\frac{n+4}2\right]+1} + \ldots + a_1^{2\left( \left[
\frac{n+4}2\right]-1\right)+1} a_{n-\left[ \frac{n+4}2\right]+1}
x_n \right)-$$

$$ - \ldots -2\beta_n' \left( a_1^{2n-2}
x_{n-1} +a_1^{2n-3} a_2x_n\right) + \beta'a_1^{2n}x_n =-2
\sum\limits_{k=\left[ \frac{n+4}2 \right]}^n
\beta_k'a_1^{2k-3}\sum\limits_{i=1}^{n-k+2} a_i x_{k+i-2} +
$$

$+ \left. 4\beta_{\left[ \frac{n+4}2\right]}'a_1^{2\left[
\frac{n-1}2\right]+1}a_{n+1}\beta_{\left[ \frac{n+4}2\right]} x_n
+ \beta' a_1^{2n} x_n. \right.$

In the case of an even $n$, for the product $\left[
y_1',y_{n+1}'\right]$ we have:
$$
\left[ y_1',y_{n+1}'\right]=-2\sum\limits_{k=\left[ \frac{n+4}2\right]}^n \beta_k' x_{k-1}' +\beta' x_n' =-2
\beta_{\left[ \frac{n+4}2\right]}' \left( a_1^{2\left(\left[ \frac{n+4}2\right]-1\right)}x_{\left[
\frac{n+4}2\right]-1}+\right.
$$ $$
+ \left. a_1^{2\left(\left[ \frac{n+4}2\right]-2\right)+1}a_2
x_{\left[ \frac{n+4}2\right]}+\ldots+ a_1^{2\left( \left[
\frac{n+4}2\right]-2\right)+1}a_{n-\left[
\frac{n+4}2\right]+1}x_n\right)- 2 \beta_{\left[
\frac{n+4}2\right]+1}' \left( a_1^{2 \left[
\frac{n+4}2\right]}x_{\left[ \frac{n+4}2\right]}+  \right.$$

$$ + \left. \left. a_1^{2\left(\left[ \frac{n+4}2\right]-1\right)+1}a_2
x_{\left[ \frac{n+4}2\right]+1}+\ldots+ a_1^{2\left( \left[
\frac{n+4}2\right]-1\right)+1}a_{n-\left[
\frac{n+4}2\right]+1}x_n\right)- \ldots - 2\beta_n' (a_1^{2n-2}
x_{n-1} + \right.$$

$$ + a_1^{2n-3} x_{n})+ \left. \beta'
a_1^{2n}x_n=-2\sum\limits_{k=\left[ \frac{n+4}2\right]}^n
\beta_k'a_1^{2k-3} \sum\limits_{i=1}^{n-k+2} a_i
x_{k+i-2}+\beta'a_1^{2n}x_n. \right.
$$

Comparing the coefficients we obtain the following restrictions:

when $n$ is odd
$$
\left\{\begin{array}{l} b_{n+1} \beta_j=a_1^{2j-3} \beta_j', \quad \left[\frac{n+4}2\right] \le j\le n, \\[2mm]
a_{n+1} b_{n+1} \gamma+a_1b_{n+1}\beta= a_1^{2n}
\beta'+4\beta_{\left[\frac{n+4}2\right]}' a_1^{2\left[
\frac{n+1}2\right]-1} a_{n+1} \beta_{\left[\frac{n+4}2\right]} \\
\end{array} \right. \eqno (3.7)
$$

when $n$ is even
$$
\left\{\begin{array}{l} b_{n+1} \beta_j=a_1^{2j-3} \beta_j', \quad \left[\frac{n+4}2\right] \le j\le n, \\[2mm]
a_{n+1} b_{n+1} \gamma+a_1b_{n+1}\beta= a_1^{2n} \beta'. \\
\end{array} \right. \eqno (3.8)
$$

It is not difficult to check that considering other
multiplications we have either restrictions (3.7)-(3.8) or
identity.

Note that from (3.6) we have $b_n=\displaystyle
\frac{-a_{n+1}b_{n+1}\gamma}{a_1}.$

Thus, combining the restrictions (3.6), (3.7) and (3.8) it follows
the proof of the theorem.
\end{proof}

Introduce the operators which are similar like $k-$dimensional
vectors:
$$
\begin{array}{rl} j & \\ V^0_{j,k}(\alpha_1, \alpha_2,\ldots, \alpha_k) =
 ( 0, \ldots, 0, 1, & \delta \sqrt[j]{\delta ^{j+1}} S_{m,j}^{j+1} \alpha_{j+1}, \delta \sqrt[j]{\delta
 ^{j+2}}
S_{m,j}^{j+2} \alpha_{j+2}, \ldots , \delta \sqrt[j]{\delta ^{k}}
S_{m,j}^{k} \alpha_{k}   ) ; \\ \end{array}
$$ $$
\begin{array}{rl} j & \\ V^1_{j,k}(\alpha_1, \alpha_2,\ldots, \alpha_k) =
 ( 0, \ldots, 0, 1, &  S_{m,j}^{j+1} \alpha_{j+1},  S_{m,j}^{j+2} \alpha_{j+2}, \ldots , S_{m,j}^{k}
\alpha_{k}  ) ; \\ \end{array}
$$ $$
\begin{array}{rl} j & \\ V^2_{j,k}(\alpha_1, \alpha_2,\ldots, \alpha_k) =
 ( 0, \ldots, 0, 1, &  S_{m,2j+1}^{2(j+1)+1} \alpha_{j+1},  S_{m,2j+1}^{2(j+2)+1} \alpha_{j+2}, \ldots ,
S_{m,2j+1}^{2k+1} \alpha_{k}  ) ; \\ \end{array}
$$ $$
V^0_{k+1,k}(\alpha_1, \alpha_2,\ldots, \alpha_k) = V^1_{k+1,k}(\alpha_1, \alpha_2,\ldots, \alpha_k) =
V^2_{k+1,k}(\alpha_1, \alpha_2,\ldots, \alpha_k) = (0, 0, \ldots, 0);
$$

$$
\begin{array}{rl} j & \\ W_{s,k}
 ( 0,  \ldots, 0, 1, & S_{m,j}^{j+1} \alpha_{j+1},
S_{m,j}^{j+2} \alpha_{j+2}, \ldots , S_{m,j}^{k} \alpha_{k}   ) =
\\ \end{array}
$$ $$
\begin{array}{rrcl} j &  & {s+j} & \\ =( 0, \ldots, 0,  1, & 0, \ldots, 0, & 1, & S_{m,j}^{j+1} \alpha_{s+j+1}, S_{m,j}^{j+2}
\alpha_{s+j+2}, \ldots, S_{m,j}^{k-s} \alpha_{k} ), \\
\end{array}
$$
$$
\begin{array}{rl} j & \\ W_{k+1-j,k}
 ( 0, \ldots,0, 1, & 0, \ldots , 0 ) =
\\ \end{array}
\\
\begin{array}{rl} j & \\ ( 0, \ldots, 0, 1, & 0, \ldots , 0 )
\\ \end{array}$$
where $k\in\mathbb C,$ $\delta=\pm 1,$ $1\le j\le k,$ $1\le s\le
k-j,$ $S_{m,t}=\displaystyle \cos\frac{2\pi m}t + i\sin \frac{2\pi
m}t $ $(m=0,1,\ldots,t-1)$.

Theorem \ref{thm32} allows us to classify the Leibniz
superalgebras from the variety $Leib^{n,m}$ with characteristic
sequence $(n \ |\ m-1, 1),$ nilindex $n+m$ and $m=n+1.$
\begin{thm} \label{thm33}
Let $L$ be a Leibniz superalgebra of variety $Leib^{n,m}$ with
characteristic sequence $(n \ |\ m-1, 1),$ nilindex $n+m$ and
$m=n+1.$ Then $L$ is isomorphic to one of the following pairwise non
isomorphic superalgebras:

if $n$ is odd (i.e. $n=2q-1$):
$$
\begin{array}{lll}
L\left(1,\delta\beta_{q+1}, V_{j,q-2}^0(\beta_{q+2}, \beta_{q+3}, \ldots, \beta_n),0\right), & \displaystyle
\beta_{q+1}\ne \pm\frac12, & 1\le j\le q-1, \\[3mm]  L\left(1,\beta_{q+1}, V_{j,q-1}^0(\beta_{q+2}, \beta_{q+3}, \ldots,
\beta_n,\beta)\right), & \beta_{q+1}=\displaystyle  \pm\frac12, & 1\le j\le q, \\[3mm] L(0,1,V_{j,q-2}^0(\beta_{q+2}, \beta_{q+3},
\ldots, \beta_n),0), & 1\le j\le q-1, & \\[3mm] L(0,0, W_{s,q-1}(V^1_{j,q-1}(\beta_{q+2}, \beta_{q+3}, \ldots, \beta_n,
\beta))), & 1\le j\le q-1, & 1\le s\le q-j, \\[2mm] L(0,0,\ldots,0); \\
\end{array}
$$

if $n$ is even (i.e. $n=2q$):
$$
\begin{array}{lll} L(1,V^2_{j,q-1}(\beta_{q+2}, \beta_{q+3}, \ldots, \beta_n,), 0), & 1\le j\le q, & \\[2mm]
L(0, W_{s,q}(V^1_{j,q}(\beta_{q+2}, \beta_{q+3}, \ldots,
\beta_n,\beta))), & 1\le j\le q, & 1\le s\le q+1-j,
\\[2mm] L(0,0,\ldots,0). \\
\end{array}
$$
\end{thm}
\begin{proof}
Consider $n$ is odd, i.e. $n=2q-1,$ where $q\in\mathbb N.$ From
Theorem \ref{thm32} we have the following restrictions:
$$
\left\{\begin{array}{l} b_{n+1}^2\gamma=\gamma'a_1^{2n}, \\
b_{n+1}\beta_j=a_1^{2j-3}\beta_j', \quad q+1\le j\le n, \\a_{n+1}
b_{n+1} \gamma+
a_1b_{n+1}\beta=a_1^{2n}\beta'+4\beta_{q+1}'a_1^na_{n+1}\beta_{q+1},
\\ \end{array} \right.
$$
for which we consider all possible cases.

\textbf{Case 1.} Let  $\gamma\ne 0.$ Then taking
$b_{n+1}=\pm\displaystyle \frac{a_1^n}{\sqrt \gamma}$, we obtain
$\gamma'=1$. Substituting the value of $b_{n+1}$ in other
restrictions we obtain equalities:
$$
\beta_{q+1+j}'=\pm\frac{\beta_{q+1+j}}{a_1^{2j} \sqrt \gamma},\quad 0\le j\le q-2,
$$ $$
\beta'= \pm\frac{a_{n+1}(\gamma-4\beta_{q+1}^2)+a_1\beta}{a_1^{n} \sqrt \gamma}.
$$

\textbf{Case 1.1.} If $\gamma-4\beta^2_{q+1}\ne 0,$ then putting $\displaystyle
a_{n+1}=-\frac{a_1\beta}{\gamma- 4\beta^2_{q+1}},$ we have $\beta'=0$ and
$\beta'_{q+1+j}=\pm\displaystyle\frac{\beta_{q+1+j}}{\sqrt \gamma}a_1^{-2j}$
for $0\le j\le q-2.$

If $\beta_{q+1+j}=0$ for any $j\in \{1, \ldots, q-2\},$ then
$\beta'_{q+1+j}=0$ and we obtain the superalgebras:
$$
L(1,\delta\beta_{q+1}, 0, \ldots, 0),\quad \delta=\pm 1.
$$

If $\beta_{q+2}=\beta_{q+3}= \ldots =\beta_{q+t}=0$ and
$\beta_{q+t+1}\ne 0$ for some $t\in \{1, 2, \ldots, q-2 \}.$ Then
taking $a_1^{-2t}= \pm\frac{\sqrt \gamma}{\beta_{q+t+1}}$ (i.e.
$a_1^{-2}= \sqrt[t]{\pm 1} \sqrt[t]{\left| \frac{\sqrt
\gamma}{\beta_{q+t+1}}\right|} \left( \cos\frac\varphi{t} +i\sin
\frac\varphi{t}\right)\left( \cos\frac{2\pi m}t+i\sin\frac{2\pi
m}t \right),$ where $\varphi= arg\left( \frac{\sqrt
\gamma}{\beta_{q+t+1}}\right),$ $m=0, 1, \ldots, t-1$), we obtain:
$$
\beta_{q+t+1}'=1\ \mbox{ and } \ \beta_{q+t+j}'= \pm \sqrt[t]{\pm
1}S_{m,t}^j \beta_{q+t+j},\quad m=0, 1, \ldots, t-1.
$$

So, in this case we have the following superalgebras:
$$
L\left(1,\delta\beta_{q+1}, V_{j,q-2}^0(\beta_{q+2}, \beta_{q+3}, \ldots, \beta_n),0\right), \quad \beta_{q+1}\ne
\pm\frac12, \quad 1\le j\le q-1,\quad \delta=\pm 1.
$$

\textbf{Case 1.2.} If  $\gamma-4\beta^2_{q+1}= 0,$ then
$\beta_{q+1}'=\displaystyle \pm\frac12$ and we have
$$
\beta_{q+1+j}'= \pm\frac{\beta_{q+1+j}}{\sqrt \gamma} a_1^{-2j},
\quad 1\le j\le q-2, \eqno (3.9)
$$ $$
\beta'= \pm\frac{\beta}{\sqrt \gamma} a_1^{-n+1}. \eqno (3.10)
$$
If we assume that $j=q-1$ in restriction (3.9), then we obtain $
\beta_{2q}'=\displaystyle \pm\frac{\beta_{2q}}{\sqrt \gamma}
a_1^{-2(q-1)}. $ Since $n=2q-1,$ then $-2(q-1)=-n+1,$ i.e. we
formally have  restriction (3.10) and therefore restriction (3.10)
can be considered as a particular case of restriction (3.9) when
$j=q-1.$

Furthermore, as in case 1.1, we obtain the following superalgebras:
$$
L\left(1,\beta_{q+1}, V_{j,q-1}^0(\beta_{q+2}, \beta_{q+3}, \ldots, \beta_n,\beta)\right), \quad
\beta_{q+1}=\displaystyle  \pm\frac12, \quad 1\le j\le q.
$$

\textbf{Case 2.} $\gamma=0.$ Then $\gamma'=0$ and
$$
b_{n+1}\beta_{q+1+j}=a_1^{2q-1+2j}\beta'_{q+1+j},\quad 0\le j\le q-2,
$$ $$
a_1b_{n+1}\beta=a_1^{2n}\beta'+4\beta'_{q+1}a^n_1a_{n+1}\beta_{q+1}.
$$

\textbf{Case 2.1.} $\beta_{q+1}\ne 0.$ Then taking $\displaystyle
b_{n+1}=\frac{a_1^n}{\beta_{q+1}}$ and $\displaystyle a_{n+1}=
\frac{b_{n+1}\beta}{4a_1^{n-1}\beta_{q+1}},$ we have
$\beta'_{q+1}=0,$ $\beta'=0$ and $\displaystyle \beta'_{q+1+j}=
\frac{\beta _{q+1+j}}{\beta_{q+1}} a_1^{-2j},$ $1\le j\le q-2.$

Furthermore,  as in case 1.1, we obtain the superalgebras:
$$
L(0,1,V_{j,q-2}^1(\beta_{q+2}, \beta_{q+3}, \ldots, \beta_n),0),
\quad 1\le j\le q-1.
$$

\textbf{Case 2.2.} $\beta_{q+1}= 0.$ Then $\beta'_{q+1}=0$ and
$$
b_{n+1}\beta_{q+1+j}=a_1^{2q-1+2j}\beta'_{q+1+j}, \quad 1\le j\le
q-2, \eqno (3.11)
$$ $$
b_{n+1}\beta=a_1^{2n-1}\beta'. \eqno (3.12)
$$

If we assume that $j=q-1,$ in restriction (3.11), then we obtain
$b_{n+1}\beta_{2q}=a_1^{2n-1}\beta'_{2q}.$ Since $n=2q- 1,$ then
$2q-1+2(q-1)=4q-3=2n-1,$ and then we formally obtain restriction
(3.12) and, therefore restriction (3.12) can be considered as a
particular case of (3.11) when $j=q-1.$

\textbf{Case 2.2.1} $\beta_{q+2}=\beta_{q+3}=\ldots=\beta_{q+t}=0$
and $\beta_{q+t+1}\ne0$ for some $t\in \{ 1,2,\ldots, q-1\}.$ Then
if we choose $b_{n+1}=\displaystyle
\frac{a_1^{2q-1+2t}}{\beta_{q+1+t}},$ we obtain $\beta_{q+1+t}'=1$
and
$$
\beta_{q+1+j}'=\frac{\beta_{q+1+j}}{\beta_{q+1+t}} a_1^{-2(j-t)}, \quad t+1\le j\le q-1.
$$

Thus, in this case we have the superalgebras:
$$
L(0,0, W_{s,q-1}(V^1_{j,q-1}(\beta_{q+2}, \beta_{q+3}, \ldots,
\beta_n, \beta))), \quad 1\le j\le q-1, \quad 1\le s\le q-j.
$$

\textbf{Case 2.2.2} $\beta_{q+j+1}=0$ for any $j$ ($1\le j\le q-1$).
Then we obtain superalgebra: $$ L(0, 0, \ldots, 0).$$

Consider the case of even $n$, i.e. $n=2q$ for some $q\in\mathbb
N.$

Then from Theorem \ref{thm32} we have the following restrictions:
$$
\left\{\begin{array}{l} b_{n+1}^2\gamma=\gamma'a_1^{2n}, \\
b_{n+1}\beta_{q+2+j}=a_1^{2q+2j+1}\beta_{q+2+j}',
\quad 0\le j\le q-2, \\a_{n+1} b_{n+1} \gamma+ a_1b_{n+1}\beta=a_1^{2n}\beta'. \\
\end{array} \right.
$$

\textbf{Case 1.}  $\gamma\ne 0.$ Then taking
$b_{n+1}=\pm\displaystyle \frac{a_1^n}{\sqrt \gamma}$ and
$a_{n+1}=\displaystyle -\frac{a_1\beta}{\gamma}$ we obtain
$\gamma'=1$; $\beta'_{q+2+j}=\pm \displaystyle
\frac{\beta_{q+2+j}}{a_1^{2j+1}\sqrt \gamma}$ $(0\le j\le q-2$) and
$\beta'=0.$

If $\beta_{q+2+j}=0$ for any $j$ ($0\le  j\le  q-2$), then we have
superalgebra:
$$ L(1,0,\ldots,0). $$

If $\beta_{q+2}=\beta_{q+3}=\ldots=\beta_{q+ t+1}=0$ and
$\beta_{q+t+2}\ne 0$ for some $t$ ($1 \le t \le q-2),$ then putting
$\displaystyle a_1^{-(2t+1)}=\pm\frac{\sqrt{\gamma}}{\beta_{q+2+t}}$
(i.e. $a_1^{-1}=\pm \sqrt[2t+1]{\left| \frac{\sqrt
\gamma}{\beta_{q+2+t}}\right|} \left( \cos\frac\varphi{2t+1} +i\sin
\frac\varphi{2t+1}\right)\left( \cos\frac{2\pi
m}{2t+1}+i\sin\frac{2\pi m}{2t+1} \right),$ where $\varphi=
arg\left(\frac{\sqrt \gamma}{\beta_{q+2+t}}\right), $ $m=0,1,\ldots,
2t),$  and substituting the value of $a_1^{-1}$ in other
restrictions we obtain
$$
\beta'_{q+2+j}=\pm\frac{\beta_{q+2+j}}{\sqrt \gamma}
\left(\pm\sqrt[2t+1]{\left|\frac{\sqrt\gamma}{\beta_{q+2+j}}
\right|}  \left( \cos\frac\varphi{2t+1} +i\sin
\frac\varphi{2t+1}\right) S_{m,2j+1}\right)^{2j+1}=
$$ $
= \beta_{q+2+j} S_{m,2j+1}^{2j+1}, \qquad t+1\le j\le q-2.
$

Thus, in this case we have the following superalgebras:
$$
L(1,V^2_{j,q-1}(\beta_{q+2}, \beta_{q+3}, \ldots, \beta_n,),0),
\quad 1\le j\le q.
$$

\textbf{Case 2.} $\gamma= 0.$ Then $\gamma'= 0$ and
$\beta'_{j+2+q}=\displaystyle \frac{b_{n+1}\beta_{q+2+
j}}{a_1^{2q+1+2j}},$ $0\le j\le q-2$, $\beta'=\displaystyle
\frac{b_{n+1}\beta}{a_1^{2n-1}}.$

Note that in this case $\beta'$ also can be considered as a
particular case of $\beta'_{j+2+q}$  for $j=q-1.$

If $\beta_{q+2+j}=0$ for any $j$ ($0 \le j \le q-1$), then we have the
superalgebra: $$ L(0, 0, 0, \ldots, 0).$$

If $\beta_{q+2}= \beta_{q+3}= \ldots= \beta_{q+ t+1}=0$ and
$\beta_{q+t+2}\ne 0$ for some $t$ ($1 \le t \le q-1$). Then putting
$\displaystyle b_{n+1}=\frac{a_1^{2q+2t+1}}{\beta_{q+2+t}},$ we
obtain:
$$
\beta'_{q+2+t}=1,\quad  \beta'_{q+2+j}=\frac{\beta_{q+2+j}}{\beta_{q+2+t}}a_1^{-2(j-t)} \quad ( t+1\le j \le q-1).
$$

As in case 2.2.1 for odd $n$, we obtain superalgebras:
$$
L(0, W_{s,q}(V^1_{j,q}(\beta_{q+2}, \beta_{q+3}, \ldots,
\beta_n,\beta))), \quad 1\le j\le q, \quad 1\le s\le q+1-j.
$$
\end{proof}

The existence of an adapted basis under the conditions of Lemma
\ref{lem31} for $m=n+2$ is represented in the following lemma.
\begin{lem} \label{lem33}
Let $L$ be a Leibniz superalgebra of variety $Leib^{n,m}$ with
characteristic sequence $(n \ |\ m-1, 1),$ nilindex $n+m$ and
$m=n+2.$ Then, there exists a basis $\{ x_1, x_2, \ldots, x_n, y_1,
y_2, \ldots, y_{n+2}\}$ of $L$ in which the multiplication has the
following form:\\[1mm]

$\begin{array}{ll}
[x_i,x_1]=x_{i+1},\ 1\le i\le n-1, & [y_j,x_1]=y_{j+1},\ 1\le j\le n, \\[1mm]
[x_i,y_1]=\frac12 y_{i+1},\ 1\le i\le n, & [y_j,y_1]=x_{j},\ 1\le j\le n,\\[1mm]
\end{array}$
$$[x_i,y_{n+2}]=\sum\limits_{k=\left[\frac{n+5}2\right]}^{n+2-i}
\beta_k y_{k-1+i},\ 1\le i\le \left[ \frac{n}2\right],  \
[y_j,y_{n+2}]=-2\sum\limits_{k=\left[\frac{n+5}2\right]}^{n+2-j}
\beta_k x_{k-2+j},\ 1\le j\le \left[ \frac{n}2\right].$$
\end{lem}
\begin{proof}
The proof of this lemma is analogous to the proof of Lemma
\ref{lem32}.
\end{proof}

Let us denote the superalgebra from the family of Lemma \ref{lem33}
by $ L\left( \beta_{\left[\frac{n+5}2\right]},
\beta_{\left[\frac{n+5}2\right]+1}, \ldots,\right.$ $\left. \beta_{n+1}\right). $

The condition of isomorphism of two superalgebras is represented in
the following theorem.
\begin{thm} \label{thm34}
Two superalgebras $L\left( \beta_{\left[\frac{n+5}2\right]},
\beta_{\left[\frac{n+5}2\right]+1}, \ldots, \beta_{n+1}\right)$
and\\ $L'\left( \beta'_{\left[\frac{n+5}2\right]},
\beta'_{\left[\frac{n+5}2\right]+1}, \ldots, \beta'_{n+1}\right)$
are isomorphic if and only if there exist $a_1,b_{n+2}\in\mathbb
C$ such that the following conditions hold:
$$
b_{n+2}\beta_j=a_1^{2j-3}\beta'_j,\quad
\left[\frac{n+5}2\right]\le j\le n+1.
$$
\end{thm}
\begin{proof}
By a change of basis the generators of the new basis are expressed
by
$$
y'_1=\sum\limits_{i=1}^{n+2}a_iy_i, \quad
y'_{n+1}=\sum\limits_{j=1}^{n+2}b_jy_j,
$$
where the rank $\left(\begin{array}{llll} a_1 & a_2 & \ldots & a_{n+2} \\ b_1 & b_2 & \ldots & b_{n+2} \\
\end{array}\right)=2,$ this alows us to express the elements of the new basis $\{
x_1',x_2',\ldots,x_n',y_1', y_2',\ldots,y_{n+2}'\}$ of the
superalgebra $L'\left( \beta'_{\left[ \frac{n+5}2\right]},
\beta'_{\left[ \frac{n+5}2\right]+1}, \ldots, \beta'_{n+1}\right)$
with respect to the elements of old basis $\{ x_1, x_2, \ldots,
x_n, y_1, $ $y_2, \ldots, y_{n+2}\}$ as:
$$
x_1'= [y_1',y_1']=a_1\sum\limits_{k=1}^n
a_kx_k-2a_{n+2}\sum\limits_{i=1}^{\left[\frac{n}2\right]} a_i
\sum\limits_{k=\left[\frac{n+5}2\right]}^{n+2-i}\beta_kx_{k-2+i};
$$
$$
x_{t+1}'=\left[ x_{t}', x_1'\right]
=a_1^{2t+1}\sum\limits_{k=1}^{n-t}a_kx_{t+k}
-2a_1^{2t}a_{n+2}\sum\limits_{i=1}^{\left[\frac{n}2\right]-t} a_i
\sum\limits_{k=\left[\frac{n+5}2\right]}^{n+2-t-i}\beta_kx_{k+t-2+i},
\ 1\le t \le \left[ \frac{n-2}2\right];
$$
$$
x_{t+1}'=\left[ x_{t}', x_1'\right] =a_1^{2t+1}\sum\limits_{k=1}^{n-t}a_kx_{t+k}, \quad \left[ \frac{n}2\right]\le
t \le n-1;
$$ $$
y_{t}'=\left[ y_{t-1}', x_1'\right] =a_1^{2(t-1)}
\sum\limits_{i=1}^{n+2-t}a_iy_{t-1+i}, \ 2\le t\le n+1.$$

If we consider the products
$$
\left[ x_{i}', x_1'\right] =\frac12 y'_{i+1}, \quad 1\le i\le n, \
\left[ y_{t}', y_1'\right] = x_t', \quad 1\le t \le n
$$
we find no restrictions.

Consider the chain of the equalities:
$$
[y_{n+2}',y_{n+2}']=b_1\sum\limits_{k=1}^n b_kx_k - 2b_{n+2}
\sum\limits_{i=1}^{\left[\frac{n}2\right]} b_i
\sum\limits_{k=\left[\frac{n+5}2\right]}^{n+2-i}\beta_kx_{k-2+i}=0.
$$

Comparing the coefficients of the basis elements in the last
equality we obtain the following restrictions:
$$
b_i=0, \ 1 \leq i \leq [\frac{n}{2}]. \eqno (3.13)
$$

From the following equalities we obtain:
$$
[y_{n+2}',y_1'] =
\left[\sum\limits_{i=[\frac{n}{2}]+1}^{n+2}b_iy_i,
\sum\limits_{j=1}^{n+2}a_jy_j\right]=a_1\sum\limits_{k=[\frac{n}{2}]+1}^{n}b_kx_k=0$$
restrictions and summing them with (3.13) we obtain $b_i=0, \ 1
\leq i \leq n.$

Therefore we have $y'_{n+2}=b_{n+1}y_{n+1}+b_{n+2}y_{n+2}.$

Consider the multiplications defining the parameters:
$$
[x_1',y_{n+2}']=\left[ a_1\sum\limits_{k=1}^n a_kx_k -2a_{n+2}\sum\limits_{i=1}^{\left[\frac{n}2\right]} a_i
\sum\limits_{k=\left[\frac{n+5}2\right]}^{n+2-i}\beta_kx_{k-2+i}, b_{n+1}y_{n+1}+b_{n+2}y_{n+2} \right]=
$$ $$
 = a_1b_{n+2}\sum\limits_{j=1}^{\left[\frac{n-2}2\right]+1} a_j \sum\limits_{k=\left[\frac{n+5}2\right]}^{n+2-j}
\beta_k y_{k+j-1}; \eqno (3.14)
$$ $$
[x_1',y_{n+2}']= \sum\limits_{k=\left[\frac{n+5}2\right]}^{n+1} \beta_k' y_k' =\beta'_{\left[\frac{n+5}2\right]}
a_1^{2\left[\frac{n+5}2\right]-2} \left( a_1 y_{\left[\frac{n+5}2\right]} +a_2 y_{\left[\frac{n+5}2\right]+1} +
\ldots + \right.
$$ $$+
\left. a_{n+2-\left[\frac{n+5}2\right]}y_{n+1}\right)+
b'_{\left[\frac{n+5}2\right]+1} a_1^{2\left[ \frac{n+5}2 \right]}
\left( a_1 y_{\left[\frac{n+5}2\right]+1} +a_2
y_{\left[\frac{n+5}2\right]+2} +\ldots + a_{n+1- \left[
\frac{n+5}2\right]}y_{n+1}\right)+
$$ $$
+ \ldots +
\beta'_{n+1}a_1^{2n+1}y_{n+1}=\sum\limits_{k=\left[\frac{n+5}2\right]}^{n+1}
\beta'_k a_1^{2(k-1)} \sum\limits_{j=1}^{n+2-k} a_j y_{k+j-1}.
\eqno (3.15)
$$

From (3.14) and (3.15) we obtain restrictions:
$$
b_{n+2}\beta_j=a_1^{2j-1}\beta'_j,\quad
\left[\frac{n+5}2\right]\le j \le n+1. \eqno (3.16)
$$
If we consider other multiplications, then we obtain restrictions
(3.16) or identities.
\end{proof}

The description up to isomorphism of a family from Lemma \ref{lem33}
is represented in the following theorem.
\begin{thm} \label{thm35}
Let $L$ be a Leibniz superalgebra of variety $Leib^{n,m}$ with
characteristic sequence $(n \ |\ m-1, 1),$ nilindex $n+m$ and
$m=n+2.$ Then $L$ is isomorphic to one of the following pairwise
non isomorphic superalgebras:
$$
L\left( W_{s,n+2-\left[\frac{n+5}2\right]} \left( V^1_{j, n+2-\left[\frac{n+5}2\right]} \left( \beta_{\left[
\frac{n+5}2 \right]}, \beta_{\left[\frac{n+5}2\right]+1}, \ldots, \beta_{n+1}\right)\right)\right),
$$
where $1\le j\le n+2-\displaystyle \left[\frac{n+5}2\right],$ $1\le s\le n+3-\displaystyle \left[
\frac{n+5}2\right]-j,$
$$ L(0,0,\ldots,0). $$
\end{thm}
\begin{proof}
From Theorem \ref{thm34} we have the following restrictions:
$$
b_{n+2}\beta'_{\left[\frac{n+5}2\right]+j}=
a_1^{2\left[\frac{n+5}2\right]+2j-3}
\beta'_{\left[\frac{n+5}2\right] +j}, \quad 0\le j\le
n+1-\left[\frac{n+5}2\right].
$$

As $\displaystyle \left[\frac{n+5}2\right]\approx q+2$, for $n=2q$ or $n=2q-1,$ then we obtain:
$$
b_{n+2}\beta_{q+j+2}=a_1^{2q+2j+1}\beta'_{q+j+2}, \quad 0\le j\le n+1-q-2.
$$

The proof of this theorem is complete by using the same arguments
as in the proof of Theorem \ref{thm33} for even case.
\end{proof}

Thus, Theorems \ref{thm33} and \ref{thm35} complete the classifications
(up to isomorphism) of Leibniz superalgebras with characteristic
sequence $C(L)=(n \ |\ m-1,1)$ and nilindex $n+m.$

{\it Acknowledgments. The last named author would like to
acknowledge the hospitality of the University of Sevilla (Spain).
He was supported by the grants INTAS - 04-83-3035 and
NATO-Reintegration ref. CBP.EAP.RIG.983169.}

{\sc Luisa M. Camacho, Jos\'{e} R. G\'{o}mez.}  Dpto. Matem\'{a}tica Aplicada I.
Universidad de Sevilla. Avda. Reina Mercedes, s/n. 41012 Sevilla.
(Spain), e-mail: \emph{lcamacho@us.es}, \emph{jrgomez@us.es}

\

{\sc Rosa M. Navarro.} Dpto. de Matem{\'a}ticas, Universidad de
Extremadura, C{\'a}ceres (Spain), e-mail: \emph{rnavarro@unex.es}

\

{\sc Bakhrom A. Omirov.} Institute of Mathematics and Information Technologues, Uzbekistan
Academy of Science, F. Hodjaev str. 29, 100125, Tashkent
(Uzbekistan), e-mail: \emph{omirovb@mail.ru}

\end{document}